\documentclass[11pt]{amsart}
\usepackage{amssymb,latexsym,amsmath,amsthm,amsfonts,amsxtra}

\theoremstyle{plain}
\newtheorem{theorem}{Theorem}

\newtheorem{lemma}{Lemma}

\theoremstyle{definition}

\theoremstyle{remark}

\numberwithin{equation}{section}

\numberwithin{equation}{section}

\newcommand{\RR}{\mathbb{R}}
\newcommand{\CC}{\mathbb{C}}

\newcommand{\B}{{\mathcal{B}}}

\newcommand{\bfn}{{\mathbf{n}}}

\def\RR{{\mathbb{R}}}
\def\CC{{\mathbb{C}}}

\def\B{{\mathcal B}}

\def\F{{\mathcal F}}

\def\H{{\mathcal H}}
\def\Id{{\hbox{\it Id}}}
\def\L{{\mathcal L}}

\def\P{{\mathcal P}}
\newcommand{\bfs}{{\mathbf{s}}}

\newcommand{\bfr}{{\mathbf{r}}}
\def\T{{\mathcal T}}
\newcommand{\bft}{{\mathbf{t}}}

\newcommand{\bfzero}{{\mathbf{0}}}

\author{Manuel Lladser}

\title[Uniqueness of polynomial canonical representations]{Uniqueness of polynomial canonical representations}

\keywords{Airy phenomena, asymptotic enumeration, analytic combinatorics, large powers of generating functions, discrete random structures, saddle point method, uniform asymptotic expansions.}

\begin{document}

\maketitle

\begin{abstract}
Let $P(z)$ and $Q(y)$ be polynomials of the same degree $k\ge1$ in the complex variables $z$ and $y$, respectively. In this extended abstract we study the non-linear functional equation $P(z)=Q(y(z))$, where $y(z)$ is restricted to be analytic in a neighborhood of $z=0$. We provide sufficient conditions to ensure that all the roots of $Q(y)$ are contained within the range of $y(z)$ as well as to have $y(z)=z$ as the unique analytic solution of the non-linear equation. Our results are motivated from uniqueness considerations of polynomial canonical representations of the phase or amplitude terms of oscillatory integrals encountered in the asymptotic analysis of the coefficients of mixed powers and multivariable generating functions via saddle-point methods. Uniqueness shall prove important for developing algorithms to determine the Taylor coefficients of the terms appearing in these representations. The uniqueness of Levinson's polynomial canonical representations of analytic functions in several variables follows as a corollary of our one-complex variables results. 
\end{abstract}

%\tableofcontents
\section{Introduction}
\label{sec:in}

Unless otherwise stated $d\ge2$ is a fixed integer and $i:=\sqrt{-1}$. We use boldface notation to denote vectors in $\CC^d$. We reserve the script $\bfzero$ to refer to the zero vector. The script $\bfr$ is reserved for a vector with strictly positive real coordinates. We refer to $\bfr$ as a polyradius. The coordinates of a vector $\bft$ are denoted $(t_1,\ldots,t_d)$. We define $\bft':=(t_1,\ldots,t_{d-1})$, in particular, $\bft=(\bft',t_d)$. The notation $|\bft|<\bfr$ means that $|t_i|<r_i$ for all $i$. Similarly, $|\bft|\le\bfr$ means that $|t_i|\le r_i$ for all $i$.\\

\noindent{\bf Problem description.} A classical example of a polynomial canonical representation is the Weierstrass preparation theorem~\cite{Tay02} which asserts the following. If $U(\bft)$ is a complex-valued analytic function in a neighborhood of $\bft=\bfzero$ and 
\[k:=\min\left\{n\ge0:\frac{\partial^n U}{\partial t_d^n}(\bfzero)\ne0\right\}<\infty\]
then there exists a polyradius $\bfr$ and analytic functions $V,u_0,\ldots,u_{k-1}$ such that
\begin{equation}
\label{ide:Weierstrass thm}
U(\bft)=V(\bft)\cdot\left\{t_d^k+\sum_{j=0}^{k-1} u_j(\bft')\cdot t_d^j\right\},
\end{equation}
for $|\bft|<\bfr$. We refer to $k$ as the {\it order of vanishing} of $U$ about the origin with respect to the variable $t_d$. The factor within the parenthesis above is called the {\it Weierstrass polynomial} of $U$ about the origin and it will be denoted as $P(\bft)$. It satisfies the following important property. For all $|\bft'|<\bfr'$ the polynomial equation in the variable $t_d$: $P(\bft',t_d)=0$, with $|t_d|<r_d$, has exactly $k$ solutions repeated according to their multiplicity. Since $V(\bfzero)\ne0$, and since the roots of a monic polynomial identify it uniquely, the factorization in (\ref{ide:Weierstrass thm}) must be unique.

The problem of whether $U(\bft)$ itself can be represented as a polynomial with respect to a possibly auxiliary variable dates back to the investigations of Chester, Friedman and Ursell~\cite{CheFriUrs57} who studied this problem for the special case of $d=2$. Later work by Levinson~\cite{Lev60a} provided a way to represent certain analytic functions of $d=2$ complex variables in a canonical way as a polynomial in two auxiliary variables. More generally, for $d\ge2$, Levinson proved the following~\cite{Lev60b}. If $U(\bft)$ is like before then there exists a polyradius $\bfr$ and analytic functions $v_0,\ldots,v_k,x$ such that
\begin{equation}
\label{ide:Levinson repr}
U(\bft)=\sum_{j=0}^k v_j(\bft')\cdot \{x(\bft)\}^j,
\end{equation}
for $|\bft|<\bfr$, with $v_j(\bfzero')=0$ for $j<k$, $v_k(\bfzero')\ne0$, and $x(\bft',0)=0$ and $\partial x/\partial t_d(\bft',0)=1$ for $|\bft'|<\bfr'$.

Unlike the representation in (\ref{ide:Weierstrass thm}), it is unclear that the representation in (\ref{ide:Levinson repr}) is unique. Indeed, the issue of uniqueness was omitted in~\cite{Lev60b} and to the best of our knowledge it has not been addressed further. The main issue surrounding the uniqueness of this representation as well as other canonical representations is the introduction of auxiliary variables. Loosely speaking, the problem is how to certify in general the validity of the following implication
\[\left[\sum_{j=0}^k v_j(\bft')\cdot x^j=\sum_{j=0}^k w_j(\bft')\cdot y^j\right]\Longrightarrow \left[v_j=w_j,\hbox{ for all $j$, and }x=y\right],\]
under the assumption that $x=x(\bft)=t_d+O(\sum_{j=1}^d|t_j|^2)$ and $y=y(\bft)=t_d+O(\sum_{j=1}^d|t_j|^2)$ for all $\bft$ sufficiently close to the origin. Clearly, to assert the uniqueness of the above factorization, it suffices to have $x=y$ in some open neighborhood of the origin in $\CC^d$. In fact, since $x(\bft',t_d)$ and $y(\bft',t_d)$ are as functions of $t_d$ locally invertible about the origin, we shall see at the end of Section \ref{sec:main results} that the validity of the above implication is closely related to the uniqueness of $y(z)=z$ as an analytic solution of the non-linear functional equation
\begin{equation}
\label{ide:fctl eq}
\left\{\begin{array}{rcl}
\prod\limits_{i=1}^k(z-z_i)&=&\prod\limits_{i=1}^k(y(z)-y(z_i)),\quad|z|<R;\\
y(0)=0&;&y'(0)=1;
\end{array}\right.
\end{equation}
where $R>0$ is a given radius and $z_1,\ldots,z_k\in\CC$ are fixed complex numbers such that $|z_i|<R$.

There is one case where the uniqueness issue of the above equation can be addressed directly. If $y(z)$ is an entire function i.e. $R=\infty$ then, according to (\ref{ide:fctl eq}), $|y(z)/z|^k\to1$ as $|z|\to\infty$. Since $y(0)=0$, $y(z)/z$ is a bounded entire function. Hence, due to Liouville's theorem~\cite{Rud87}, $y(z)/z$ must be constant and therefore $y(z)=z$ because $y'(0)=1$. Unfortunately, the case with $R=\infty$ is not of much use to address uniqueness issues of polynomial canonical representations because --- almost always --- they only apply locally. \\

%THIS IS NOT TRUE IF k=2: Furthermore, under special circumstances, the equation in (\ref{ide:fctl eq}) may have non-trivial solutions when $R<\infty$.

\noindent{\bf Connections with mixed powers generating functions.} Polynomial canonical representations are pivotal for analyzing the asymptotic behavior of oscillatory integrals~\cite{BleHan86}. Integrals of this type arise frequently in the context of asymptotic enumeration or the analysis of discrete random structures~\cite{PemWil05}.

A mixed power generating function is a generating function of the form $\prod_{i=1}^d\{f_i(z)\}^{n_i}$, where the factors $f_1,\ldots,f_d$ are complex-valued analytic functions near $z=0$ and $n_1,\ldots,n_d$ are nonnegative integers. The term of mixed power was introduced in~\cite{Lla06b} to emphasize the fact that one is usually interested in the coefficient of $z^{n_0}$ of $\prod_{i=1}^d\{f_i(z)\}^{n_i}$ as $\max\{n_0,n_1,\ldots,n_d\}\to\infty$. If one defines $\bfn:=(n_1,\ldots,n_d)$, this is equivalent to request that $\|(n_0,\bfn)\|\to\infty$ where $\|\cdot\|$ is any norm in $\RR^{1+d}$.

Generating functions of the above form occur naturally in the context of the Lagrange inversion formula~\cite{GouJac04,Wil90} with $d=1$. More recent applications include the case of $d=2,3$ to analyze the core size of various types of random planar maps~\cite{BanFlaSchSor01}.

Coefficients of mixed powers generating functions have been considered in the literature for factors $f_i$ with nonnegative coefficients by Drmota~\cite{Drm94a}, for $d=1$ and $n_0,n_1\to\infty$ at a comparable rate. Gardy~\cite{Gar95} considered the case of nonnegative coefficients for $d\ge1$ with $n_0=\Theta(n_1)$ or $n_0=o(n_1)$ and $n_i=o(n_1/\sqrt{n_0})$ for $i>1$. A geometrically based approach, in the lines used by Pemantle and Wilson~\cite{PemWil02,PemWil04}, was proposed in~\cite{Lla06b} to handle factors $f_i$ with possibly negative Taylor coefficients. Given $(t_0,\bft)\in\RR^{1+d}$ with nonnegative coordinates and such that $\|(t_0,\bft)\|=1$, say that $x$ is a {\it strictly minimal critical point} associated with $(t_0,\bft)$ provided that
\begin{eqnarray*}
t_0&=&\sum_{i=1}^d t_i\cdot\frac{xf_i'(x)}{f_i(x)};\\
\prod_{i=1}^d|f_i(z)|^{t_i}&<&\prod_{i=1}^d|f_i(x)|^{t_i},
\end{eqnarray*}
for all $z$ such that $|z|=|x|$ and $z\ne x$. If the above conditions hold and some pathological behavior is ruled out, it follows from~\cite{Lla06b} that
\begin{equation}
\label{ide:saddle int}
[z^{n_0}]\prod_{i=1}^d\{f_i(z)\}^{n_i}\sim\frac{x^{-n_0}}{2\pi}\prod_{i=1}^d\{f_i(x)\}^{n_i}\cdot\int_{-\pi}^\pi\exp\Big\{-\|(n_0,\bfn)\|\cdot F(\theta;(t_0,\bft))\Big\}d\theta,
\end{equation}
for $(n_0,\bfn)$ such that $(n_0,\bfn)/\|(n_0,\bfn)\|=(t_0,\bft)$, as $\|(n_0,\bfn)\|\to\infty$. The function $F$ is a computable function that is continuous in its $(d+2)$ arguments however it is also analytic in the variable $\theta$. For a fixed $(t_0,\bft)$, it satisfies that $F=\partial F/\partial\theta=0$ at $\theta=0$, and the real-part of $F$ is minimized at $\theta=0$. Furthermore, the above expansion applies uniformly for all $(n_0,\bfn)/\|(n_0,\bfn)\|\in\mathbb{T}$, provided that $\mathbb{T}$ is a compact set such that for all $(t_0,\bft)\in\mathbb{T}$, $x=x(t_0,\bft)$ is a strictly critical point associated with $(t_0,\bft)$ that depends continuously on $(t_0,\bft)$. In particular, the asymptotic analysis of the above integral is amenable for the saddle-point method to obtain uniform asymptotic expansions for the coefficients in (\ref{ide:saddle int}) for $(n_0,\bfn)\in\|(n_0,\bfn)\|\cdot\mathbb{T}$, as $\|(n_0,\bfn)\|\to\infty$.

It is precisely for the asymptotic analysis of integrals such as the one occurring in~(\ref{ide:saddle int}) that polynomial canonical representations of the type in (\ref{ide:Levinson repr}) play a crucial role. In particular, uniqueness of these representations shall prove important to determine the Taylor coefficients of the various terms and auxiliary variables occurring in these representations. This should aid 
in automatizing the extraction of asymptotic formulae for coefficients of mixed powers generating functions as well as multivariable generating functions. 

The lack of analyticity of $F$ in (\ref{ide:saddle int}) with respect to the variable $(t_0,\bft)$ can be over passed by thinking of $F$ as a function of $(\theta;(t_0,\bft);x)$. The original function $F$ can then be recovered by evaluating this new function at $(\theta;\bft;x(t_0,\bft))$. In order to apply the saddle-point method let $k$ be the order of vanishing of $F$ about $(0;(t_0,\bft);x(t_0,\bft))$ with respect to the variable $\theta$. Since $F=\partial F/\partial\theta=0$ at points of this type, it follows from~\cite{Lla06a} that Levinson's polynomial cannonical representation takes the form
\[F(\theta;(s_0,\bfs);x)=\sum_{j=2}^k v_j((s_0,\bfs);x)\cdot \{y(\theta;(s_0,\bfs);x)\}^j,\]
with $y=0$ and $\partial y/\partial \theta=1$ at points of the form $(0;(s_0,\bfs);x)$ that are near $(0;(t_0,\bft);x(t_0,\bft))$. If $k=2$ the above translates into having the integral appearing in (\ref{ide:saddle int}) to be described asymptotically by the Gamma function. In particular, the integral is of order $\|(n_0,\bfn)\|^{-1/2}$ as $\|(n_0,\bfn)\|\to\infty$. On the other hand, if $k=3$ the integral is described asymptotically by the Airy function. In this case the integral in (\ref{ide:saddle int}) has typically an asymptotic series expansion which is a linear combination of terms of order $(\|(n_0,\bfn)\|^{-l-1/3})_{l\ge0}$ and also of order $(\|(n_0,\bfn)\|^{-l-2/3})_{l\ge0}$. See~\cite{BleHan86,Lla03} to follow up on uniform asymptotics for integrals that involve the Gamma and Airy function.

The interested reader is referred to~\cite{Lla06b} for concrete applications of the above methodology with $k=2,3$. The reader is also referred to~\cite{BanFlaSchSor01} for a related yet more specialized discussion with $k=3$. Although our motivation to study the uniqueness of polynomial canonical representations has been argued in the context of mixed powers generating functions, they also play a fundamental role in the extraction of asymptotics of multivariable generating functions. The reader is referred to~\cite{PemWil02,PemWil04,Lla06a} to follow up on this last remark.

%SECTION
\section{Main results}
\label{sec:main results}

We first introduce two one-complex variable results. Theorem~\ref{thm:A} provides sufficient conditions to ensure that all the roots of a polynomial $Q(y)$ are contained in the range of an analytic function $y(z)$ when there exists another polynomial $P(z)$, of the same degree as $Q(y)$, such that $P(z)=Q(y(z))$ in a neighborhood of $z=0$. Under an appropriate rescaling, this translates into having $\prod_{i=1}^k(z-z_i)=\prod_{i=1}^k(y(z)-y(z_i))$, where $k$ is the degree of $P(z)$ and $z_1,\ldots,z_k$ are the roots of $P(z)$ repeated according to their multiplicity. Theorem~\ref{thm:B} provides sufficient conditions in order to conclude from this that $y(z)=z$. Both theorems are then used to show the uniqueness of Levinson's representation in (\ref{ide:Levinson repr}). The proofs of our main two theorems are presented in Section~\ref{sec:proofs}. Our main results and accompanying proofs are refined versions of some of the results obtained by the author in his doctoral dissertation~\cite{Lla03}.\\

%SUBSECTION
\noindent{\bf Auxiliary results.} In what follows, $R>0$ is a given radius and we use the notation
\begin{eqnarray*}
D&:=&\{z\in\CC:|z|<R\},\\
\H&:=&\{y:D\to\CC\hbox{ such that $y$ is analytic}\}.
\end{eqnarray*}
For $0\le r<R$, we define
\[\|f\|_r:=\sup_{|z|\le r}|f(z)|=\sup_{|z|=r}|f(z)|,\]
where the last identity is justified by the maximum modulus principle~\cite{Rud87}.

%THEOREM
\begin{theorem}
\label{thm:A}
Let $P$ and $Q$ be polynomials of the same degree $k\ge1$ and assume that $D$ contains all the roots $z_1,\ldots,z_k$ of $P$ repeated according to their multiplicity. If $y\in\H$ is such that
$P(z)=Q(y(z))$, for all $z\in D$, then 
\begin{equation}
\label{ide:equiv fact}
\left[Q^{-1}\{0\}\subset y(D)\right]\Longleftrightarrow\left[\forall i:y'(z_i)\ne0,\hbox{ and }\forall i,j:y(z_i)=y(z_j)\Leftrightarrow z_i=z_j\right].
\end{equation}
Furthermore, if either of the conditions in (\ref{ide:equiv fact}) apply then there exists a constant $q\in\CC$ such that
\begin{equation}
\label{ide:fact Q}
Q(y)=q\cdot\prod_{i=1}^k(y-y(z_i)).
\end{equation}
\end{theorem}

%THEOREM
\begin{theorem}
\label{thm:B}
For all $\rho$ and $r$ such that $0\le 2\rho<r<R$ there exists a $\delta>0$ such that if $\max_i|z_i|\le\rho$ then $y(z)=z$ is the only solution of the non-linear functional equation
\begin{equation}
\label{thm:ide:fctl eq}
\left\{\begin{array}{rcl}
\prod\limits_{i=1}^k(z-z_i)&=&\prod\limits_{i=1}^k(y(z)-y(z_i)),\quad y\in\H;\\
y(0)&=&0,
\end{array}\right.
\end{equation}
that satisfies the condition $\|y(z)-z\|_r\le\delta$.
\end{theorem}

%SUBSECTION
\noindent{\bf Proof of uniqueness of Levinson's representation.} We use the stated theorems to show the uniqueness of Levinson's polynomial canonical representation in (\ref{ide:Levinson repr}). Thus consider $U(\bft)$ analytic in a neighborhood $\bft=\bfzero$ and assume that
\begin{equation}
\label{ide:two Levinson}
U(\bft)=\sum_{j=0}^k v_j(\bft')\cdot s^j=\sum_{j=0}^k w_j(\bft')\cdot t^j
\end{equation}
where $v_j(\bfzero')=w_j(\bfzero')=0$ for $j<k$, $v_k(\bfzero')\ne0$, $w_k(\bfzero')\ne0$, and $s=t=0$ and $\partial s/\partial \theta_d=\partial t/\partial \theta_d=1$ at all points in the domain of $s$ and $t$ of the form $(\bft',0)$. We show that $v_j=w_j$, for all $j$, and that $s=t$. For this consider the transformation $\Phi(\bft)=(\bft',s(\bft))$. Since $\Phi(\bfzero)=\bfzero$ and the Jacobian matrix of $\Phi$ at $\bfzero$ is lower-triangular with non-zero entries along the diagonal, the inverse mapping theorem~\cite{Tay02} implies that $\Phi^{-1}$ is a well-defined analytic function in some open neighborhood of the origin in $\CC^d$. Define $V(\bft',z):=U(\Phi^{-1}(\bft',z))$ and $x=x(\bft',z):=t(\Phi^{-1}(\bft',z))$. It follows from (\ref{ide:two Levinson}) that
\begin{equation}
\label{ide:zx Levinson}
V(\bft',z)=\sum_{j=0}^k v_j(\bft')\cdot z^j=\sum_{j=0}^k w_j(\bft')\cdot x^j.
\end{equation}
Observe that $x=0$ and $\partial x/\partial z=1$ at all points in the domain of $x$ of the form $(\bft',0)$. Furthermore, according to the first identity above, $V$ vanishes to degree $k$ about the origin in the variable $z$. In particular, the Weierstrass preparation theorem~\cite{Tay02} implies that, for all $\bft'$ sufficiently close to $\bfzero'$, the roots of $V(\bft',z)$ can be listed as $z_1(\bft'),\ldots,z_k(\bft')$, repeated according to their multiplicity. Since for $\bft'$ sufficiently close to the origin the transformation $z\to x(\bft',z)$ is a one-to-one transformation, we may use Theorem~\ref{thm:A} in (\ref{ide:zx Levinson}) to conclude that
\[v_k(\bft')\cdot\prod_{j=1}^k\{z-z_j(\bft')\}=w_k(\bft')\cdot\prod_{j=1}^k\{x(\bft',z)-x(\bft',z_j(\bft'))\}.\]
But observe that, according to (\ref{ide:zx Levinson}), $x(\bfzero',z)=z\cdot(v_k(\bfzero')/w_k(\bfzero'))^{1/k}$ provided that the appropriate branch for the $k$-th root is selected. With this choice of branch, introduce the auxiliary variable
\[y=y(\bft',z):=x(\bft',z)\cdot\left\{\frac{v_k(\bft')}{w_k(\bft')}\right\}^{-1/k}.\]
Notice that
\[\prod_{j=1}^k\{z-z_j(\bft')\}=\prod_{j=1}^k\{y(\bft',z)-y(\bft',z_j(\bft'))\},\]
for all $\bft'$ sufficiently close to the origin in $\CC^{d-1}$ and $z$ such that $|z|<R$, where $R>0$ is certain real parameter independent of $\bft'$. But observe that, according to the Weierstrass preparation theorem, if $\bft'$ is sufficiently close to the origin then $|z_j(\bft')|<R/4$, for all $j$. On the other hand, since $y(\bfzero',z)=z$ and $y$ is uniformly continuos over compact subsets of its domain, it follows for $r=3R/4$ that
\[\lim_{\bft'\to\bfzero'}\left\|y(\bft',z)-z\right\|_r=0.\]
Theorem~\ref{thm:B} implies that $y(\bft',z)=z$, for all $\bft'$ sufficiently close to the origin and all $z$ such that $|z|<R$. In particular, $x(\bft',z)=z\cdot(v_k(\bft')/w_k(\bft'))^{1/k}$. Since $\partial x/\partial z=1$ at all points in the domain of $x$ of the form $(\bft',0)$, we conclude that $x(\bft',z)=z$. This in (\ref{ide:zx Levinson}) implies that $v_j=w_j$, for all $j$. Furthermore, since $x(\bft',z):=t(\Phi^{-1}(\bft',z))$, with $\Phi(\bft)=(\bft',s(\bft))$, we find $s=t$. This shows that Levinson's polynomial canonical representations are unique.\qed

%SECTION
\section{Proofs of main results}
\label{sec:proofs}

%SUBSECTION
\noindent{\bf Proof of Theorem~\ref{thm:A}.} Assume that $Q^{-1}\{0\}\subset y(D)$ i.e. that all roots of $Q(y)$ lie in the range of the function $y(z)$. Then the roots of $Q(y)$ may be listed as $y(\xi_1),\ldots,y(\xi_k)$ --- repeated according to their multiplicity --- in such a way that $\xi_i=\xi_j$ if and only if $y(\xi_i)=y(\xi_j)$. Define $y_i:=y(\xi_i)$ and let $n_i$ be the multiplicity of $y_i$ as a root of $Q(y)$. Observe that
\begin{equation}
\label{ide:multiplicities}
\lim_{z\to\xi_i}\frac{P(z)}{(z-\xi_i)^{n_i}}=\lim_{z\to\xi_i}\left\{\frac{y(z)-y(\xi_i)}{z-\xi_i}\right\}^{n_i}\cdot\frac{Q(y(z))}{(y(z)-y_i)^{n_i}}=\{y'(\xi_i)\}^{n_i}\cdot\lim_{y\to y_i}\frac{Q(y)}{(y-y_i)^{n_i}},
\end{equation}
where for the last identity we have used that $y_i$ is in the interior of $y(z)$ as asserted by the open mapping theorem~\cite{Rud87}. Since the limit on the right-hand side above exists, $\xi_i$ has to be a root of $P(z)$ of multiplicity at least $n_i$. In particular, using that the multiplicity of $\xi_i$ in the list $\xi_1,\ldots,\xi_k$ is precisely $n_i$, $\prod_{i=1}^k(z-\xi_i)$ must divide $P(z)$. Since $P(z)$ is of degree $k$, $\prod_{i=1}^k(z-\xi_i)=\prod_{i=1}^k(z-z_i)$ and as a result the sequence $\xi_1,\ldots,\xi_k$ is just a reordering of $z_1,\ldots,z_k$. From this it is immediate that $y(z_i)=y(z_j)$ if and only if $z_i=z_j$, and that $Q(y)$ factorizes as described in (\ref{ide:fact Q}). Furthermore, since $\xi_i$ must be a root of multiplicity $n_i$ of $P(z)$, it follows from (\ref{ide:multiplicities}) that  $y'(\xi_i)\ne0$. In particular, $y'(z_i)\ne0$ for all $i$.

To complete the proof of the theorem it suffices to show that if for all $i$ and $j$, $y'(z_i)\ne0$, and $y(z_i)=y(z_j)$ if and only if $z_i=z_j$, then all roots of $Q(y)$ lie in the range of $y(z)$. For each $i$ let $m_i$ be the multiplicity of $z_i$ as a root of $P(z)$ and define $y_i:=y(z_i)$. Since, like in (\ref{ide:multiplicities}), we have
\[\lim_{y\to y_i}\frac{Q(y)}{(y-y_i)^{m_i}}=\left\{\frac{1}{y'(z_i)}\right\}^{m_i}\cdot\lim_{z\to z_i}\frac{P(z)}{(z-z_i)^{m_i}}\]
and the limit on the right-hand side above is nonzero and finite, it follows that $y_i$ is a zero of multiplicity $m_i$ of $Q(y)$. Since the multiplicity of $y_i$ in the list $y_1,\ldots,y_k$ is $m_i$, $\prod_{i=1}^k(y-y(z_i))$ divides $Q(y)$. In particular, all roots of $Q(y)$ lie in the range of $y(z)$ because $Q(y)$ is of degree $k$. This completes the proof of the theorem.\qed\\

%SUBSECTION
\noindent{\bf Proof of Theorem~\ref{thm:B}.} To prove this result it would suffice to show that the differential of the operator associated with (\ref{thm:ide:fctl eq}) is one-to-one and has a continuous inverse. However, a technical issue with this approach is that it is unclear that the inverse of the differential is continuous at all when considering the natural Banach space $\B:=\{y\in\H:y\hbox{ extends continuously to the boundary of $D$}\}$ with the infinite-norm $\|y\|_\infty:=\sup_{|z|\le R}|y(z)|$. Technical difficulties arise even showing that the pre-images of functions in $\B$ under the differential stay bounded near the boundary of $D$.

Due to the above considerations we consider a weaker topology. We embed $\H$ with the topology of uniform convergence over compact subsets of $D$. In particular, a sequence $y_1,y_2,\ldots\in\H$ converges to $y\in\H$ provided that $\lim_{n\to\infty}\|y-y_n\|_r=0$, for all $0\le r<R$. This topology is induced by a metric in $\H$ under which this space is complete i.e. Cauchy sequences are convergent~\cite{Tay02}.
We define $\H_0:=\{y\in\H:y(0)=0\}$. Clearly $\H_0$ is a closed linear subspace of $\H$. In particular, $\H_0$ is also complete when endowed with the topology of uniform convergence over subsets of $D$.

%LEMMA
\begin{lemma}\label{lem:isomorphism}
If $2\cdot\max_i|z_i|<R$ then the operator $\L:\H_0\to\H$ defined as
\[(\L f)(z):=\frac{1}{k}\sum\limits_{j=1}^k \frac{f(z)-f(z_j)}{z-z_j}\]
is a linear isomorphism.
\end{lemma}

\noindent{\it Proof.} Define $\rho:=\max_i|z_i|$. According to the removable singularity theorem~\cite{Tay02}, $\L$ is a well-defined linear transformation. Since on the other hand, for $f\in\H_0$ and $\rho<r<R$ it applies that
\[\|\L f\|_r\le2\|f\|_r\cdot\sup_{|z|=r}\frac{1}{k}\sum_{j=1}^k\left|\frac{1}{z-z_j}\right|\le\frac{2\|f\|_r}{r-\rho},\]
it is immediate that $\L$ is a continuous linear operator.

To show that $\L$ is one-to-one consider the polynomial $p(z):=\prod\limits_{j=1}^k(z-z_j)$ and observe that
\begin{equation}
\label{ide:p'}
p'(z)=p(z)\cdot\sum_{j=1}^k\frac{1}{z-z_j}.
\end{equation}
Suppose that $f\in\H_0$ is such that $\L f=0$. Using (\ref{ide:p'}), a simple calculation reveals that
\[f(z)\cdot p'(z)=\sum_{j=1}^kf(z_j)\cdot \frac{p(z)}{z-z_j}.\]
Since $p'(z)$ is a polynomial of degree $(k-1)$ in the variable $z$ and the right-hand side above is a polynomial of degree at most $(k-1)$, the division algorithm implies that there is a constant $c\in\CC$ and a polynomial $q(z)$ of degree at most $(k-2)$ such that $f(z)=c+q(z)/p'(z)$, for $|z|<R$ such that $p'(z)\ne0$.
But it is well-known that the roots of $p'(z)$ are convex linear combinations of $z_1,\ldots,z_k$. In particular, the $(k-1)$ roots of $p'(z)$ lie in the disk of radius $\rho$ centered at the origin. Since $f(z)$ is bounded in this disk however $q(z)$ is of degree at most $(k-2)$, it follows from the above identity that $q=0$. Hence $f$ is constant and therefore $f=0$ because $f\in\H_0$. This shows that $\L$ is one-to-one.

To show that $\L:\H_0\to\H$ is an isomorphism we define an operator $\T:\H\to\H_0$ such that $\L(\T f)=f$, for all $f\in\H$. With the understanding that $0^0=1$ define
$\alpha_n:=\sum_{j=1}^kz_j^n/k$. Observe that $|\alpha_n|\le\rho^n$. In particular, the series $A(z)=\sum_{n=0}^\infty\alpha_nz^n$ is analytic for $|z|<1/\rho$. Since $A(0)=1$, $1/A(z)$ is analytic in a neighborhood of the origin. Hence if $\beta_n$ is defined as the coefficient of $z^n$ of the power series expansion of $1/A(z)$ about the origin, it applies that
\begin{equation}
\label{ide:convolution}
\sum\limits_{j=i}^n\beta_{n-j}\cdot\alpha_{j-i}=\left\{
\begin{array}{lcl}
0 &,& n>i\,,\\
1 &,& n=i\,.
\end{array}\right.
\end{equation}
Furthermore, an inductive argument shows that
\begin{equation}
\label{ide:bound betas}
|\beta_n|\le (2\rho)^n.
\end{equation}
For $f\in\H$ with $f(z)=\sum_{n=0}^\infty f_nz^n$ define
\begin{equation}
\label{ide:def T}
(\T f)(z):=z\cdot\sum_{j=0}^\infty\left\{\sum_{l=j}^\infty f_l\cdot\beta_{l-j}\right\}z^j.
\end{equation}
To see that the above transformation is well-defined consider $2\rho<r_0<r_1<R$. According to the Cauchy estimates~\cite{Rud87}, $|f_n|\le\|f\|_{r_1}\cdot r_1^{-n}$. As a result, using (\ref{ide:bound betas}) we obtain that
\begin{equation}
\label{ide:bound coeff T}
\left|\sum_{l=j}^\infty f_l\cdot\beta_{l-j}\right|\le\frac{\|f\|_{r_1}}{1-2\rho/r_1}\cdot r_1^{-j}.
\end{equation}
Since the above holds for any $2\rho<r_1<R$, $(\T f)(z)$ is analytic for $|z|<R$ and $(\T f)(0)=0$. Furthermore, $\T$ is a continuous linear operator because the above inequality implies that
\[\|\T f\|_{r_0}\le\frac{r_0}{(1-2\rho/r_1)(1-r_0/r_1)}\cdot\|f\|_{r_1}.\]
Finally we show that $\L(\T f)=f$, for all $f\in\H$. For this let $\P$ be the linear subspace of polynomials in the complex variable $z$. Since $\T$ is continuous and $\P\subset\H$ is a dense linear subspace, it suffices to show that $\L(\T f)=f$, for all $f\in\P$. However, due to linearity, this is equivalent to show that $\L(\T z^n)=z^n$, for all $n\ge0$. For this observe that
\[\L(\T z^n)=\L\left(\sum_{j=0}^n\beta_{n-j}\,z^{j+1}\right)=\sum_{j=0}^n\beta_{n-j}\,\L(z^{j+1})=\sum_{j=0}^n\beta_{n-j}\sum_{i=0}^j\alpha_{j-i}\,z^i=\sum_{i=0}^n\left\{\sum_{j=i}^n\beta_{n-j}\alpha_{j-i}\right\}z^i=z^n,\]
where for the last equality we have used (\ref{ide:convolution}). This completes the proof of the lemma. \qed

%LEMMA
\begin{lemma}
\label{lem:inequality}
Let $0\le 2\rho<r<R$. If $\max_i|z_i|\le\rho$ then $\|(\L f)\|_r\ge c\cdot\|f\|_r$, for all $f\in\H_0$, with
\begin{equation}
\label{ide:constant c}
c:=\min_{(z_1,\ldots,z_k)}\inf_{|z|=r}\left|\frac{1}{k}\sum_{j=1}^k\frac{1}{z-z_j}\right|\cdot\left\{1+\frac{\rho/r}{(1-2\rho/r)^3}\right\}^{-1}>0,
\end{equation}
where the minimum is taken over all $(z_1,\ldots,z_k)$ such that $\max_i|z_i|\le\rho$.
\end{lemma}

\noindent{\it Proof.} Let $\T:\H\to\H_0$ be the inverse operator of $\L$. According to  (\ref{ide:p'}), since $\L(\T f)=f$ for all $f\in\H$, it applies for $|z|<R$ that
\[(\T f)(z)\cdot p'(z)=p(z)\cdot\left\{k\cdot f(z)+\sum_{j=1}^k\frac{(\T f)(z_j)}{z-z_j}\right\}.\]
Define $c_1:=\sup_{|z|=r}|p(z)/p'(z)|$ and observe that $0<c_1<\infty$ because $r>\rho$. It follows from the above identity that
\[\|\T f\|_r\le c_1\cdot k\cdot\left\{\|f\|_r+\frac{\|\T f\|_\rho}{r-\rho}\right\}
\le  c_1\cdot k\cdot\left\{1+\frac{\rho}{(1-2\rho/r)^3}\right\}\cdot\|f\|_r,\]
where for the second inequality we have used (\ref{ide:def T}) and (\ref{ide:bound coeff T}) with $r_1=r$. The above implies that for all $f\in\H_0$, $c_2\cdot\|f\|_r\le\|\L f\|_r$, with
\[0<c_2:=\frac{1}{k\cdot c_1}\left\{1+\frac{\rho/r}{(1-2\rho/r)^3}\right\}^{-1}=\inf_{|z|=r}\left|\frac{1}{k}\sum_{j=1}^k\frac{1}{z-z_j}\right|\cdot\left\{1+\frac{\rho/r}{(1-2\rho/r)^3}\right\}^{-1},\]
where for the second identity we have used (\ref{ide:p'}). Since for $z_1,\ldots,z_k,y_1,\ldots,y_k$ such that $\max_i|z_i|\le\rho$ and $\max_i|y_i|\le\rho$ it applies that
\[\left|\inf_{|z|=r}\left|\frac{1}{k}\sum_{j=1}^k\frac{1}{z-z_j}\right|-\inf_{|z|=r}\left|\frac{1}{k}\sum_{j=1}^k\frac{1}{z-y_j}\right|\right|\le\frac{\max_j|z_j-y_j|}{(r-\rho)^2},\]
$c_2$ depends continuously on $(z_1,\ldots,z_k)$. This shows (\ref{ide:constant c}) and completes the proof of the lemma.\qed\\

We finally prove Theorem~\ref{thm:B}. Let $0\le2\rho<r<R$ be such that $\max_i|z_i|\le\rho$. To study the uniqueness of the functional equation in (\ref{thm:ide:fctl eq}) consider the operator $\F:\H_0\to\H$ defined as $(\F y)(z):=\prod_{i=1}^k\{y(z)-y(z_i)\}$. Given $y\in\H_0$ define $f:=y-Id$. Observe that
\[\F(y)-\F(\Id)=\F(f+\Id)-\F(\Id)=p(z)\cdot\sum_{J}\prod_{j\in J}\frac{f(z)-f(z_j)}{z-z_j},\]
where the index $J$ in the summation varies over all possible non-empty subsets of the set $\{1,\ldots,k\}$. In particular,
\[\F(y)-\F(\Id)=p(z)\cdot\left\{\sum_{j=1}^k\frac{f(z)-f(z_j)}{z-z_j}+\sum_{I}\prod_{i\in I}\frac{f(z)-f(z_i)}{z-z_i}\right\},\]
where the index $I$ varies over all possible subsets of the set $\{1,\ldots,k\}$ of cardinality two or greater. As a result, we find that
\begin{equation}
\label{ide:main inequality}
\|\F(y)-\F(\Id)\|_{r}\ge\inf_{|z|=r}|p(z)|\cdot\left\{\left\|\sum_{j=1}^k\frac{f(z)-f(z_j)}{z-z_j}\right\|_r-\left\|\sum_{I}\prod_{i\in I}\frac{f(z)-f(z_i)}{z-z_i}\right\|_r\right\}.
\end{equation}
According to Lemma~\ref{lem:inequality}, we have that
\[\left\|\sum_{j=1}^k\frac{f(z)-f(z_j)}{z-z_j}\right\|_r\ge k\cdot c\cdot\|f\|_r,\]
for an appropriate constant $c$ which depends on $\rho$ and $r$ but not on $z_1,\ldots,z_k$. On the other hand, if $\|f\|_r\le1$ we also have that
\[\left\|\sum_{I}\prod_{i\in I}\frac{f(z)-f(z_i)}{z-z_i}\right\|_r\le\sum_{i=2}^k{k\choose i}\frac{2^i\|f\|_r^i}{(r-\rho)^i}\le\left(1+\frac{2}{\rho}\right)^k\|f\|_r^2.\]
As a result, if $0<\|f\|_r\le\min\left\{1,k\cdot c\cdot(1+2/\rho)^{-k}\right\}$, it follows from (\ref{ide:main inequality}) that $\|\F(y)-\F(\Id)\|_{r}>0$. This completes the proof of Theorem~\ref{thm:B}.\qed

%\subsection{Proper Names}
%\label{sec:names}
%\acknowledgements
%\label{sec:ack}

%\nocite{*}
%\bibliographystyle{abbrvnat}
% use the following instead if you encounter problems 
\bibliographystyle{alpha}
\bibliography{biblio}
%\label{sec:biblio}

\end{document}